\documentclass[12pt]{article}

\usepackage[margin=2.5cm]{geometry}

\usepackage{amsmath, amssymb, latexsym, amsthm}

\usepackage{graphicx}
\usepackage[dvipsnames]{color}

\newtheorem{theorem}{Theorem}[section]
\newtheorem{proposition}[theorem]{Proposition} 

\newtheorem{corollary}[theorem]{Corollary}
\newtheorem{lemma}[theorem]{Lemma}

\newtheorem*{theorem*}{Theorem}
\newtheorem*{conjecture*}{Conjecture}
\newtheorem*{corollary*}{Corollary}
\newtheorem*{proposition*}{Proposition}

\newtheorem{problem}[theorem]{Problem}
\newtheorem{claim}[theorem]{Claim}

\theoremstyle{definition}

\newtheorem{example}{Example}

\theoremstyle{remark}

\newcommand{\thetaint}{\theta_{\mathrm{int}}}

\begin{document}

\title{Decomposing graphs into interval colorable subgraphs
and no-wait multi-stage schedules}
 %with applications in scheduling}

\author{Armen S. Asratian, Carl Johan Casselgren, Petros A. Petrosyan}

\author{
{\sl Armen S. Asratian}\thanks{{\it E-mail address:} 
armen.asratian@liu.se}\\ 
Department of Mathematics \\
Link\"oping University \\ 
SE-581 83 Link\"oping, Sweden
\and
{\sl Carl Johan Casselgren}\thanks{{\it E-mail address:} 
carl.johan.casselgren@liu.se}\\ 
Department of Mathematics \\
Link\"oping University \\ 
SE-581 83 Link\"oping, Sweden
\and
{\sl Petros A. Petrosyan}\thanks{{\it E-mail address:} 
petros\_petrosyan@ysu.am} \\
Department of Informatics \\ and Applied Mathematics,\\
Yerevan State University \\ 0025, Armenia
}

\maketitle

\begin{abstract}
A graph $G$ is called \emph{interval colorable} if it has a proper edge coloring with colors 
$1,2,3,\dots$ such that the colors of the edges incident to every
vertex  of $G$ form an interval of integers. Not all  graphs are interval colorable;
in fact, quite few families have been proved to admit interval colorings.
%;for instance, odd cycles do not admit interval colorings.
In this paper we introduce and
investigate a new notion,  the interval coloring thickness
of a graph $G$, denoted $\thetaint(G)$, which is  the
minimum number of  interval colorable edge-disjoint subgraphs of 
$G$ whose union is $G$. 

Our investigation is motivated by scheduling problems 
with compactness requirements,
in particular, problems whose solution may consist of several schedules,
but where each schedule must not contain any waiting periods or idle times for
all involved parties.
 We first prove that every connected   properly $3$-edge colorable graph  with
maximum degree $3$ is interval
colorable, and using this result,
we deduce an upper bound on  $\theta_{int}(G)$
for general graphs $G$. 
We demonstrate that this upper bound can be improved in the case when $G$ is bipartite, 
 planar 
or complete multipartite  and consider
some applications in  timetabling.
\end{abstract}

	\section{Introduction}
	
	\subsection{Background, motivation and our contribution}
	
	The classical graph coloring problem is the problem of assigning positive integers,
	identified as colors, to the vertices or edges of a graph
	so that no two adjacent vertices/edges
	receive the same color. A variety of topics in
	computer science and operations research such as
	scheduling, frequency assignment, and register
	allocation may be formulated as graph coloring problems, see e.g.
	\cite{Briggs, Gandham, Marx} and references therein.
	Many concrete applications involve extra constraints.
	Consider, for example, the classic school timetabling problem
	with an additional compactness requirement.

	\begin{example}
	In a school we want to schedule lectures so that they are consecutive
	for both teachers and classes.
	Constructing a schedule with $t$ time periods satisfying these requirements  is
	equivalent to the problem of finding an  edge coloring of the graph $G$,
	with vertices for teachers and classes, and where every edge 
	represents a lecture given
	by a certain teacher to a certain class,
	with  colors  $1,2,\dots,t$   such that the colors of the edges incident 
	to every vertex
	of $G$  are distinct and form an interval of integers.
	Such a coloring is 
	called an {\em interval $t$-coloring} of $G$. 
	\end{example}

	The notion of interval colorings was introduced by
	Asratian and Kamalian in 1987 \cite{AsrKam}
	(available in English as \cite{AsrKamJCTB}),
	motivated by the problem of constructing such timetables.
	A  graph is {\em  interval colorable} if it has an interval 
	$t$-coloring for some integer $t$.
	Not all graphs are interval colorable,
	a simple example is the complete graph $K_3$ with $3$ vertices.
	
	Now consider another scenario.
	
\begin{example}
\label{ex:matches}
Suppose that a provincial children's
soccer competition is to take place in a town during
a couple of days. In such a competition, teams from all over the province participate
and each team takes part in a limited number of matches against competing teams.
Moreover, since the competition involves teams from all over the province,
for practical reasons,
we require that every team's matches
should be consecutive without any waiting periods during each of the competition days.
Is it possible to schedule all matches during $k \geq 1$ days so that the
matches are consecutive for every team during each day?
\end{example}
	
In general, a solution to such a scheduling problem
as in Example	\ref{ex:matches},
where we ask for a schedule in a total of $t$ time units partitioned
into $k \geq 1$ {\em stages}, i.e. disjoint time periods, 
so that a ``no-wait'' condition
holds at each stage, we call a
{\em no-wait multi-stage schedule}.

We can model Example \ref{ex:matches} in graph theoretical terms
by forming a graph where vertices represent teams and edges represent matches.
(Note that this graph is not complete, since not every team plays against every other team).
Then a no-wait multi-stage schedule with $k$ stages exists if and only if
there is a {\em decomposition} of $G$ into $k$ interval colorable subgraphs,
that is, a list of $k$ subgraphs of $G$ such that every edge of $G$ appears in exactly
one subgraph in the list.
Thus, in general terms, the question of the existence of a no-wait multi-stage schedule
can be formulated in graph theoretical terms as follows.
	
\begin{problem}
\label{prob:decomp} Let $G$ be a graph and $k$ a positive integer.
	Is there a decomposition of $G$ into $k$ interval 
	colorable subgraphs?
\end{problem}

The minimum integer $k$ for which $G$ admits such a decomposition we call
{\em the interval coloring thickness} of $G$ and denote by $\thetaint(G)$.

In this paper, we introduce and investigate the parameter $\thetaint(G)$;
our investigation is motivated by scheduling problems
whose solutions permit a schedule partitioned into several stages, but each 
such stage must not contain any waiting periods or idle times for each involved party.
	
For instance, in the timetabling problem described above, we might ask
for a {\em weekly} school timetable, where the schedule for each day
satisfies the no-wait condition; this application is investigated in some
detail in Section 
\ref{timetable}. Naturally, there are many more
scheduling problems where a no-wait multi-stage schedule is desirable;
let us here just consider two further concrete examples.	
Firstly, we have the following variation of the well-known open shop problem 
\cite{Gonzalez}
with an additional compactness requirement.
	
\begin{example}
	We are given
	$m$ processors $P_1,\dots,P_m$ and $n$ jobs $J_1,\dots,J_n$ to 
	be processed within a period of $t$ consecutive time units. 
	  Each job $J_i$ consists  of  tasks $T_{i1},\dots,T_{im}$ 
		which has to be processed
	  on $P_1,\dots,P_m$, respectively. 
		For simplicity, we assume that
		the processing time  of each task is $0$ or $1$.
	   Different tasks of the same job cannot be processed simultaneously 
		and no processor can work on two tasks at the same time.
	Furthermore, we assume that  all tasks $T_{i1},\dots,T_{im}$ of $J_i$ 
	are to be executed contiguously, and, similarly, there should be
	no waiting periods for the processors, i.e.
	all tasks of $P_j$ should be executed contiguously.
	
	Is it possible to schedule all jobs
	within a total of $t$ consecutive time units, partitioned into
	$k$ stages, so that  
	the no-wait condition holds at each stage? 
	That is, rather than
	requiring that the whole schedule satisfies a  no-wait condition, which is
	a quite strong requirement, we
	ask for a no-wait multi-stage schedule. Naturally, such a schedule is desirable
	when a production process may be partitioned into 
	several disjoint time periods, e.g.
	days.
	A solution with $k$ stages exists if and only if the
	bipartite graph $B$ with vertices $J_1, \dots, J_n$,
	$P_1,\dots, P_m$, where the vertices $J_i$ and $P_j$ are
	joined by an edge if the processing time  of the task $T_{ij}$ is $1$,
	satisfies $\thetaint(B) \leq k$. 
	\end{example}

	  %An example of Problem \ref{problem2} is the problem of existence  of  a weekly 
		%school timetable  for $k$ days 
%without interruptions for teachers and classes,  considered in section \ref{timetable}.

	%Our investigation in this paper is motivated by scheduling problems 
	%whose solutions consist of several schedules where each schedule 
	%satisfies a ``no-wait'' condition.

Finally, let us consider the following variation of a problem described 
by Bodur and Luedtke \cite{BodurLuedtke}.

\begin{example} 
%In \cite{BodurLuedtke}, Bodur and Luedtke considered the job interviews model with compactness requirements. Here we consider some variation of that model. 
Suppose that some firms organise job interviews for possible candidates during a couple of days. We need to provide the schedule of job interviews where neither firm representatives nor candidates wait between their meetings during these days. If we construct a bipartite graph $H$ with parts $F$ and $C$, where vertices in $F$ represent firms and vertices in $C$ represent candidates, and edges represent the required interviews, then the minimum number of days needed for a schedule of job interviews without waiting periods is precisely equal to 
 $\thetaint(H)$. 
\end{example}

From a theoretical point of view, the  problem of 
determining whether a (bipartite) graph 
$G$ has an interval coloring (or, equivalently, whether $\thetaint(G)=1$) is NP-complete \cite{Seva}.
%so the problem of determining $\thetaint(G)$ for a general graph $G$ is NP-complete.
	%Generally, it is an NP-complete problem
	%to determine whether a (bipartite) graph
	%has an interval coloring \cite{Seva}. 
	However some classes of %bipartite
	graphs have been proved to admit interval colorings; it is known,
	for example, that trees, regular and complete bipartite graphs
	\cite{AsrKam,Hansen},
	bipartite graphs with maximum degree at most three \cite{Hansen},
	doubly convex bipartite graphs
	\cite{AsrDenHag,KamDiss}, grids \cite{GiaroKubale1}, 
	outerplanar bipartite graphs \cite{GiaroKubale2} 
	and some classes of biregular bipartite graphs 
	\cite{West, CarlJToft, 
	CasselgrenPetrosyanToft,
	Hansen, HansonLotenToft, KamMir, Pyatkin, YangLi}
	have interval colorings.

Due to the NP-hardness of computing the parameter  $\thetaint(G)$, 
in this paper we focus on locating
tractable instances and proving constructive upper bounds on $\thetaint$
for different families of graphs.
%to determine ``Is $\thetaint(G)=1$?'' (or, equivalently,  
%``Is $\thetaint(G)>1?$'') is    NP-complete even for bipartite graphs.
First we prove that every connected 
$3$-edge colorable
graph with maximum degree $3$ is interval colorable, and
%It was shown in \cite{AsrKamJCTB} that the condition $\chi'(G)=\Delta(G)$ is a necessary (but not sufficient) condition for a graph $G$ to be interval colorable. Here we show that this condition is  sufficient for $G$ 
%to be interval colorable if $G$ is a  connected graph  with $\Delta(G)\leq 3$. 
using this result we 
deduce a general upper bound on the interval coloring thickness of an arbitrary  graph.
Then we investigate the parameter $\thetaint$ for various families of graphs;
in particular, we
demonstrate how this general upper bound
on $\thetaint(G)$ can be improved in the case when $G$ is bipartite, planar
or complete $r$-partite. 
%For the case when the input graph is bipartite, we also
%consider a specific scheduling problem in some detail.
We conclude the paper by  pointing to some
open problems.

\subsection{Graph theoretical preliminaries}

	We use \cite{West1} for terminology and notation not defined here. 
	All graphs considered are
  finite, undirected, allow multiple edges and contain no loops, unless otherwise stated.

A \emph{simple} graph is a graph
without loops and multiple edges.
$V(G)$ and $E(G)$ denote the sets of vertices and
edges of a graph $G$, respectively.
We denote by $\Delta(G)$ and $\delta(G)$ the maximum and minimum  degrees of
the vertices of a graph $G$, respectively, and by $d_G(v)$ (or just $d(v)$)
the degree of a
vertex $v$ in $G$. A graph $G$ is called {\em subcubic} if $\Delta(G)\leq 3$.
The {\em distance $d_G(u,v)$} between two vertices $u$ and $v$ in $G$ is the length
of a shortest path between $u$ and $v$ in $G$.

A \emph{$2$-factor} of a graph $G$ (where loops are allowed) is a
$2$-regular spanning subgraph of $G$. Petersen's well-known
theorem \cite{Petersen}
asserts that every $2r$-regular
graph (where loops are allowed) can be decomposed
into edge-disjoint $2$-factors.

		A proper $t$-edge coloring of a graph $G$ is a
    mapping $\alpha:E(G)\longrightarrow  \{1,\dots,t\}$ such that
    $\alpha(e)\not=\alpha(e')$ for every pair of adjacent
    edges $e$ and $e'$ in $G$. If $e\in E(G)$ and $\alpha(e)=k$ then
    we say that the edge $e$ is {\em colored $k$}.
    The {\em chromatic index $\chi'(G)$} of a graph $G$
    is the minimum number $t$ for which there exists a proper
    $t$-edge coloring of $G$.	
		The two famous theorems on edge coloring are K\"onig's theorem \cite{Konig}, 
		which states
		that $\chi'(G) = \Delta(G)$ if $G$ is bipartite, and
		Vizing's theorem \cite{Vizing} which asserts that 
		$\chi'(G)\leq \Delta(G)+1$ for any simple graph $G$.
		In addition, Shannon's theorem \cite{Shannon} states that 
		$\chi'(G) \leq \frac{3}{2} \Delta(G)$
		for any graph $G$.
		Vizing's theorem partitions the set of simple graphs into two sets, namely the
		graphs $G$ that satisfy $\chi'(G) = \Delta(G)$, 
		and $\chi'(G) = \Delta(G)+1$, respectively;
		the latter family of graphs is called {\em Class 2} and the former 
		{\em Class 1}.
		
		If $\alpha $ is an edge coloring of $G$ and $v\in V(G)$, then
$S_{G}\left(v,\alpha \right)$ (or $S\left(v,\alpha \right)$) denotes
the set of colors appearing on edges incident to $v$.
For two positive integers $a$ and $b$ with $a\leq b$, we denote by
$\left[a,b\right]$ the interval of integers $\left\{a,a+1,\ldots,b\right\}$.

\section{A new class of interval colorable graphs}

It was shown in \cite{AsrKamJCTB} that the condition $\chi'(G)=\Delta(G)$ 
is  necessary (but not sufficient) for a graph $G$ to be interval colorable. 
Here we show that this condition is  sufficient for $G$ to be interval colorable if
$G$ is a connected subcubic graph; 
this generalizes the result of Hansen
\cite{Hansen} that all bipartite subcubic graphs are interval colorable.

\begin{theorem}
\label{th:general}
	Let $G$ be a graph with $\chi'(G)=\Delta(G) \leq 3$.
	If no component of $G$ is an odd cycle, then $G$ is interval colorable.
\end{theorem}

Clearly, the condition on odd cycles is necessary.
Before proving
the theorem, let us also note that it is sharp 
with respect to the maximum degree;
the graph obtained from a triangle $v_1v_2v_3v_1$ by adding
a path of length two with a new internal vertex vertex $v_{ij}$
between any pair of distinct vertices $v_i, v_j$ satisfies the condition $\chi'(G)=\Delta(G)=4$, but is not interval colorable.

\begin{proof}[Proof of Theorem \ref{th:general}]
	If $G$ has maximum degree at most $2$, 
	then $G$ is trivially
	interval colorable,  so assume that $\Delta(G)=3$.
	 Since $\chi'(G)=\Delta(G)$,  there is a proper $3$-edge coloring of $G$
	using colors $1,2,3$.
	Let $M$ be the set of edges colored $1$.
	
	From this edge coloring of $G$, we define an edge-colored graph $T$
	as follows: a vertex $v$ is in $V(T)$ if and only if $v \in V(G)$
	and $v$ is a vertex of degree $1$ in $G-M$, or $v$ is incident
	with an edge of $M$ (or both).
	The edge set of $T$ consists of blue, red and green edges:
	\begin{itemize}
		
		\item every edge of $M$ is in $E(T)$ and is colored red;
		
		\item two vertices of degree $1$ that are endpoints of a (maximal) path in
		$G-M$ of even length are joined by an blue edge in $T$;
		
		\item two vertices of degree $1$ that are endpoints of a (maximal) path in
		$G-M$ of odd length are joined by an green edge in $T$.
	
	\end{itemize}
	
	Note that   $T$ is a bipartite graph with maximum degree $2$.
	Moreover, if $C$ is a cycle of $T$, then every second edge of $T$ is red.
	
	\begin{claim}
	\label{cl:vertexcol}
		There is a vertex coloring $c :V(T) \to \{A,B\}$ such that
		every edge that is in a path or in a cycle of $T$ with an even number of blue edges,
		satisfies that
		\begin{itemize}
		
		\item[(i)] both endpoints of a red or green edge are colored by the same color
		($A$ or $B$);
		
		\item[(ii)] the endpoints of a blue edge are colored by different colors.
		
		\end{itemize}
		Moreover, every edge  of a cycle $C$ with an odd number
		of blue edges satisfies (i) and (ii) except for one red edge $e$ of $C$
		whose endpoints are colored by different colors. Additionally, we can
		select this edge $e$ to be any of the red edges of $C$.
		\end{claim}
	\begin{proof}
		The claim trivially holds for any path or cycle with an even number of
		blue edges; just color the vertices of every
		component of $T$ sequentially so that they satisfy the conclusion of the claim.
		That the statement holds for cycles of $T$ with an odd number of blue edges is
		a straightforward exercise left to the reader.
	\end{proof}
	
	We continue the proof of the theorem.
	Let $c$ be a coloring of $T$ satisfying the conclusion of the claim.
	
	Consider a cycle $C$ in $T$ with an odd number of blue edges.
	If $C$ corresponds to an even cycle in $G$, then $C$ is {\em good}.
	If $C$ corresponds to an odd cycle $C_G$ in $G$, then by assumption,
	some edge of $M$ must be incident with a vertex of $C_G$; that is,
	some vertex of $C_G$ of degree two in $G-M$
	is in $T$. Consider a maximal path $P$ in $G-M$ that
	is in $C_G$ and such that some internal vertex $u$ of $P$ is incident
	with an edge of $M$. We choose the coloring $c$ on $C$
	such that a red edge
	incident with an endpoint $v$ of $P$ is colored by different colors.
	Moreover, without loss of generality, we assume that $u$ is
	an internal vertex of $P$ incident with an edge of $M$, such that
	the distance between $v$ and $u$ is minimal. Finally, assume
	that if $d_G(u,v)$ is odd, then $c(u) \neq c(v)$, and if
	$d_G(u,v)$ is even, then $c(u)=c(v)$. (This can be achieved by possibly
	swapping all colors (under $c$) 
	on the component $C$ if necessary; note that $u$ is not in $C$.)

	Thus, we may assume that for every cycle $C$ in $T$ with an odd number of blue
	edges, if $C$ is not good, then the following holds:
	
	\begin{itemize}
	
	\item[(a)] there is exactly one red edge $e$ of $C$ whose endpoints are colored 
	differently. Moreover, one endpoint of $e$ is an endpoint of
	a component $P$ in $G-M$ that is a path, which satisfies
	that some internal vertex of $P$ is incident with an edge from $M$. 	
	Additionally, if $u$ is the internal vertex of $P$ with shortest distance to $v$ that is
	incident with an edge from $M$, then $u$ and $v$ are colored 
	differently (under $c$)
	if $d_P(u,v)$ is odd, and they are colored by the same color (under $c$)
	if $d_P(u,v)$ is even.
	\end{itemize}
	
	We are now ready to color the edges of $G$; since $G-M$ is a bipartite
	graph with maximum degree two, its components are paths and even cycles.
	
	For every component $Q$ of $G-M$, we color the edges alternately by colors
	$2$ and $3$. If $Q$ is a cycle then this yields a proper coloring of the
	edges of $Q$; if $Q$ is a path
	then we start with color $2$ at a vertex colored $A$ if $Q$, and color
	$3$ at a vertex of $Q$ colored $B$. 
	This yields a proper coloring of $G-M$ since every odd maximal path of
	$G-M$ corresponds to a green edge in $T$, and every even maximal path
	of $G-M$ corresponds to a blue edge.
	
	An edge of $M$ is colored by color $1$ if both its endpoints are colored $A$,
	and by color $4$ if both its endpoints are colored $B$.
	Since an endpoint $x$ of a path $P$ in $G-M$ is incident to an edge colored $2$
	from $P$ if $c(x)=A$, and an edge colored $3$ if $c(x)=B$,
	this yields an interval coloring of the hitherto edge-colored subgraph of $G$.
	
	It remains to color the edges of $M$ that correspond to red edges 
	whose endpoints are colored differently under $c$.
	Every such red edge
	is in a unique cycle in $T$ with an odd
	number of blue edges. So consider such a cycle $C$ in $T$.
	If $C$ is good, then we recolor all its edges by colors $2$ and $3$ alternately.
	If $C$ is not good, then we consider some cases.
	
	Let $e=vx$ be an edge of $M$ in $C$, whose endpoints are colored differently.
	Then $v$ is an endpoint of a path $P$
		in $G-M$ containing an internal vertex $u$, that is incident with an edge from $M$,
		and satisfying that $d_P(u,v)$ is minimal w.r.t to this property.
	
	\begin{itemize}
		
		\item 
		If $d_P(u,v)$ is even and $c(u)=c(v)=A$, then
		we recolor the edges on the portion of $P$ from $u$ to $v$ 
		by colors $0$ and $1$ alternately, and starting with color
		$0$ at $u$, and color $xv$ by the color $2$;
		if $c(u)=c(v)=B$, then we proceed similarly, but use colors $5$ and $4$
		on the portion of $P$ from $u$ to $v$, and starting with color $5$ at $u$,
		and color $xv$ by the color $3$;
		
		\item if $d_P(u,v)$ is odd and $c(u)=A$ and $c(v)=B$, then
		we recolor the edges on the portion of $P$ from $u$ to $v$ 
		by colors $0$ and $1$ alternately, and starting with color
		$0$ at $u$, and color $xv$ by the color $1$;
		if $c(u)=B$ and $c(v)=A$, then we proceed similarly, but use colors $5$ and $4$
		on the portion of $P$ from $u$ to $v$, and starting with color $5$ at $u$,
		and color $xv$ by the color $4$.
	
	\end{itemize}

	Recall that red edges of $T$  whose endpoints are colored differently under $c$ 
	are all in different
	components of $T$. Hence,
	by repeating this process for every red edge of $T$ whose endpoints are colored
	differently (under $c$) we obtain an interval edge coloring of $G$.
	\end{proof}

\begin{corollary}
A connected subcubic graph $G$ is interval colorable if and only if $\chi'(G)=\Delta(G)$.
\end{corollary}

By K\"onig's edge coloring theorem \cite{Konig}, $\chi'(G) = \Delta(G)$ for any 
 bipartite graph $G$. This and Theorem \ref{th:general} imply the following:

\begin{corollary}(\cite{Hansen})
\label{Hansen}
	If $G$ is a bipartite graph with maximum degree at most $3$,
	then $G$ has an interval coloring.
\end{corollary}

Other families of graphs that satisfy
the hypothesis of Theorem \ref{th:general} include the so-called series-parallell 
graphs, that is, graphs that  can be obtained 
recursively from a single edge by the operations of subdividing and doubling edges.
Every simple series-parallel graph with maximum degree at least 3 is Class 1 
\cite{Nishizeki}; 
thus every such graph
is interval colorable. 
Since every outerplanar graphs is series-parallell, 
this generalizes a result of \cite{ArsPetros}.

\section{Upper bounds on the interval coloring thickness of graphs}
	
A natural strategy for proving upper bounds on the interval
coloring thickness of graphs is to consider decompositions into graphs
that are known to admit interval colorings.
In this section we present some results in this vein.
	
	\subsection{A general upper bound on $\thetaint(G)$}

	In this section we prove a general upper bound on $\theta_{int}(G)$ 
	for an arbitrary graph $G$ in terms of its chromatic index.
	For the proof of this result we need the following lemma, 
	the proof of which is left to the reader.

\begin{lemma}
\label{lem:add}
	Let $H$ be an interval colorable graph. If $G$ is obtained by 
	\begin{itemize}
	\item[(i)] adding a new  pendant edge to $H$, or 
	
	\item[(ii)] adding a cycle $C$ with exactly
	one common vertex $v$ with $H$ that has degree $1$ in $H$, 

	\end{itemize}
	then
	$G$ is interval colorable.
\end{lemma}

  Note that Lemma \ref{lem:add} implies the following:

\begin{proposition}
\label{prop:cactus}
	If $G$ is a connected  graph with $\Delta(G)\geq 3$ where any two cycles are vertex-disjoint, then
	$G$ admits an interval coloring.
\end{proposition}

Let us now prove the main result of this section.
	
\begin{theorem}
  \label{th:five}
		 For any graph $G$,
		 $\thetaint(G) \leq 2\left\lceil \frac{\chi'(G)}{5} \right\rceil$.
  \end{theorem}
  
  \begin{proof}
  Without loss of generality we assume that $G$ is connected. 
  
  {\bf Case 1}. $\chi'(G)\leq 5$.
  
  We will show that $\theta_{int}(G)\leq 2$.
  Consider a proper 5-edge coloring of $G$ and let $M_i$ denote the set 
	of edges colored $i$, for $i=1,\dots,5$.
  Furthermore, let $H$ and $F$ be the subgraphs of $G$ induced by the 
	sets $M_1\cup M_2$ and $M_3\cup M_4\cup M_5$, respectively.
  Clearly, $\Delta(F)\leq 3$ and $H$ is a bipartite graph with 
	$\Delta(H)\leq 2$.
	
	Let $F_0$ be the subgraph of $F$ obtained from $F$ by
	removing all components from $F$ that are odd cycles.
	Then $F_0$ is interval colorable by Theorem \ref{th:general}.

	Now, if $F$ contains no odd cycles, then we are done,
	because $H$ is interval colorable.
	Otherwise, the nontrivial components of
	$F -E(F_0)$ are $l\geq 1$ uncolored components $C_i$
	that are odd cycles.

 Since $G=H\cup F$ is connected, there is a 
 path $P_1$ in $H$ where one endpoint $x_1$ 
is in $V(F_0)$, and the other endpoint $y_1$ is 
in one of the odd cycles $C_1,\dots, C_l$, say $C_1$. 
Moreover, we assume that $C_1$ is the odd cycle 
from $\{C_1, \dots, C_l\}$
with shortest
distance to $F_0$ in $G$; then $P_1$ is disjoint
from all the cycles $C_2,\dots, C_l$.

By Lemma \ref{lem:add}, the graph 
$F_0\cup C_{1} \cup P_1$ is interval colorable.
Put $F_1=F_0\cup C_{1} \cup P_1$ and $H_1:=H-E(P_1)$.
Then $H_1$ is interval colorable as well.
The number of components that are odd cycles in 
$F-V(F_1)$ is $l-1$; thus
by continuing this process we will obtain an interval 
colorable graph $F_l$ containing
all odd cycles of $F$.
Moreover, the subgraph $H_l$ of $H$ 
is interval colorable because it is
bipartite and has maximum degree at most $2$.
Hence, $\thetaint(G) \leq 2$.
\bigskip

{\bf Case 2}. $\chi'(G)\geq 6$.

Consider a proper edge coloring of $G$ with colors 
$1,2,\dots,\chi'(G)$.
Let $t=\left\lfloor \frac{\chi'(G)}{5}\right\rfloor$ and let $M_i$ denote the set of edges colored $i$,  for 
$i=1,\dots,\chi'(G)$. 
We define the subgraphs $G_1,\dots,G_k$ of $G$, where $k=\left\lceil \frac{\chi'(G)}{5} \right\rceil$ as follows:   
For each $j=1,\dots,t$,  let $G_j$ be the subgraph   induced by the set of edges  
$M_{5j-4}\cup M_{5j-3}\cup\dots \cup M_{5j}$. 
If $t=k$, then $G = \bigcup_{i=1}^t G_i$; otherwise
if $t\not=k$, then $t+1=k$. 
In this case we define $G_{k}$ as the subgraph  induced by the set 
$\cup_{i=t+1}^k M_i$.
For each  $i=1,\dots,k$, $G_i$ is a graph with $\Delta(G_i)\leq 5$ colored with 5 colors. Therefore, 
by the result in Case 1, % Theorem \ref{th:five}, 
$\thetaint(G_i)\leq 2$. Thus $\thetaint(G)\leq 2k=2\left\lceil \frac{\chi'(G)}{5} \right\rceil.$
\end{proof}

By the well-known edge coloring theorems by Shannon \cite{Shannon} and Vizing \cite{Vizing} we have the following
consequence of	
Theorem \ref{th:five}. 
 
\begin{corollary}
\label{bound}

(i) For any graph $G$,  $\thetaint(G)\leq 2\left\lceil \frac{3\Delta(G)}{10}\right\rceil$.

 (ii) For any simple Class 1 graph $G$, $\thetaint(G) \leq 2\left\lceil \frac{\Delta(G)}{5} \right\rceil$.
 
 (iii) For any simple Class 2 graph $G$, $\thetaint(G) \leq 2\left\lceil \frac{\Delta(G)+1}{5} \right\rceil$. 
\end{corollary} 

It is well-known that the upper bounds in Vizing's and Shannon's theorems can be achieved
by polynomial-time algorithms (see e.g. \cite{ScheideStiebitz}), 
so the proof of Theorem \ref{th:five} combined with
these proofs yield efficient algorithms for finding decompositions and colorings
attaining the upper bounds in Corollary \ref{bound}.
\bigskip

\subsection{Bipartite graphs}
	
	There is a prominent line of research on interval colorings of bipartite graphs.
	This family of graphs is also particularly interesting due to applications
	in scheduling problems and timetabling. In this section we prove upper bounds
	on $\thetaint(G)$ for bipartite graphs $G$
	 and 
	also consider a specific application involving weekly school timetables.
	\subsubsection{General upper bounds on $\thetaint(G)$} 

Every bipartite graph $G$ with $\Delta(G)\leq 3$ is, by Corollary \ref{Hansen}, 
interval colorable. It is also easy to show that
 $\thetaint(G) \leq \left\lceil \frac{\Delta(G)}{3} \right\rceil$ 
for any bipartite graph $G$ with $\Delta(G)\geq 4$ by proceeding 
as in the proof of
Theorem \ref{th:five}:
		 partition $E(G)$ into $\Delta=\Delta(G)$ matchings $M_1,\dots,M_{\Delta}$ and
		define $k=\left\lceil \frac{\Delta}{3} \right\rceil$ disjoint 
		sets of edges consisting of
		  triples of matchings, and possibly one consisting of one or two matchings.
		  Each such set of edges induces a 
			bipartite subgraph of maximum degree not exceeding 3, which thus is interval
			colorable by Corollary \ref{Hansen}.

Therefore $\theta_{int}(G)\leq k$.
		   
		We  prove here a slightly stronger result which can be useful in scheduling
		problems where several schedules of roughly equal ``size''
		is desirable	(see Proposition \ref{prop:even}).

	\begin{proposition}
	\label{th:Bip}
		If $G$ is a bipartite graph with $\Delta(G)\geq 4$, 
		then $\thetaint(G) \leq \left\lceil \frac{\Delta(G)}{3} \right\rceil$. Moreover 
		$G$ can be decomposed into
		$\left\lceil \frac{\Delta(G)}{3} \right\rceil$ edge-disjoint interval colorable subgraphs
		in such a way  that at each vertex the numbers of incident edges, 
		in any pair of subgraphs, differ by at most one.
		\end{proposition}
		
	\begin{proof}
			Let %$G$ be a bipartite graph and
			 $k=\left\lceil \frac{\Delta(G)}{3} \right\rceil$.
			We define a new bipartite graph $H$ by splitting each vertex $v$  of degree 
			at least $k+1$ into $\left\lfloor \frac{d(v)}{k} \right\rfloor$ vertices of 
			degree $k$ and possibly 
			(if $d(v)>k\cdot  \left\lfloor \frac{d(v)}{k} \right\rfloor$)
			 one extra vertex of degree 
			$d(v)-k\cdot  \left\lfloor \frac{d(v)}{k} \right\rfloor$.
			 The partitioning of the edges in this splitting is arbitrary, 
			other than ensuring that each vertex receives the correct degree.
			 Since the graph $H$ has maximum degree $k$, 
			it has, by K\"onig's edge coloring theorem \cite{Konig},  
			a proper $k$-edge-coloring. 
			Now collapse $H$ back into $G$ and consider the induced  
			coloring of the edges of $G$.
			Since each vertex of $G$ is split into at most $3$
			 vertices, every $v\in V(G)$ can be incident to at most 
			$\left\lceil \frac{d(v)}{k} \right\rceil$
			edges of each color.	On the other hand, corresponding to each vertex $v$
			there were $\left\lfloor \frac{d(v)}{k} \right\rfloor$ 
			vertices of degree $k$.
			Thus there must be at least 
			$\left\lfloor \frac{d(v)}{k} \right\rfloor$ edges of 
			each color incident to $v$. 
			This means that the obtained coloring of $G$ 
			satisfies the following condition:  
			at each vertex the numbers of incident edges, in any pair of colors, 
			differ by at most 1.
			Let $G_i$ be the subgraph of $G$ induced by the edges of color $i$, 
			$1\leq i\leq k$.
			By construction,  the degree of each vertex in $G_i$ does not exceed 3.
			Then, by 		
			Corollary \ref{Hansen},  $G_i$ is interval colorable, for $i=1,2,\dots,k$. 
			Therefore $ \theta_{int}(G) \leq 
			\left\lceil \frac{\Delta(G)}{3} \right\rceil$.
			\end{proof}

	We remark that since there are well-known efficient algorithms 
	for finding an optimal 
	proper edge coloring of a bipartite graph, the proof of the preceding 
	proposition
	combined with the proof of Theorem \ref{th:general},
	yields
	a polynomial-time algorithm for finding a decomposition, 
	and interval colorings of the corresponding parts, attaining the upper bound
	 $\thetaint(G) \leq \left\lceil \frac{\Delta(G)}{3} \right\rceil$.
	This also holds for our next result, which is an improvement of the 
	preceding proposition
	for Eulerian bipartite graphs.
	To prove it we need the following lemma.

\begin{lemma}
\label{Lem:add} If $G$ is a bipartite graph with
$\Delta(G)=2r$ and the degree of every vertex in $G$ is $1,2$ or $2r$, then $G$ has an 
interval $2r$-coloring $\alpha$ such that for each 
$v\in V(G)$ with $d_{G}(v)=2$, $S_{G}(v,\alpha)=\{2i-1,2i\}$ for some $i$, $1 \leq i \leq r$.
\end{lemma}

The proof of this lemma is similar to the proof of Theorem 2.1
in \cite{AsratianCasselgrenPetrosyan}, so we shall omit it.

\begin{theorem}
\label{th:Eul}
If $G$ is an Eulerian bipartite graph, then $\theta_{int}(G)\leq \left\lceil \frac{\Delta(G)}{4} \right\rceil.$	
\end{theorem}

\begin{proof}[Proof of Theorem \ref{th:Eul}]
For the proof, we construct a new graph $G^{\star}$ as follows: for each vertex $u \in V(G)$ of degree $2k$, we add $\frac{\Delta(G)}{2}-k$ loops at $u$ 
	$\left(1\leq k< \frac{\Delta(G)}{2}\right)$. Clearly, $G^{\star}$ is a $\Delta(G)$-regular
graph. By  Petersen's theorem \cite{Petersen}, $G^{\star}$ can be decomposed into
a union of edge-disjoint $2$-factors $F_{1},\ldots,F_{\frac{\Delta(G)}{2}}$. By removing all
loops from the $2$-factors $F_{1},\ldots,F_{\frac{\Delta(G)}{2}}$, we
obtain that the resulting graph $G$ is a union of edge-disjoint
subgraphs $F^{\prime}_{1},\ldots,F^{\prime}_{\frac{\Delta(G)}{2}}$, where each
$F'_{i}$ is a collection of even cycles in $G$.

Since the maximum degree in $G$ is even, $\Delta(G) = 4l-2$ or $\Delta(G)=4l$, for some
integer $l$.
We now define the subgraphs $Q_1, \dots, Q_l$, by setting $Q_i = F'_{2i-1} \cup F'_{2i}$,
$i=1,\dots,l-1$, and setting $Q_l= F'_{2l-1} \cup F'_{2l}$ if $\Delta(G) = 4l$,
and $Q_l = F'_{2l-1}$ if $\Delta(G) =4l-2$. By Lemma
\ref{Lem:add}, each
$Q_{i}$  has an interval coloring; thus $\theta_{int}(G)\leq l$.

\end{proof}

		For a bipartite graph $G$ with  parts  $X$ and $Y$, we denote
		the maximum degree of the vertices in $X$ by $\Delta(X)$;
		$\delta(X)$ denotes the minimum degree of these vertices.
		
		For bipartite graphs where one of the parts $X$ and $Y$ 
		has small maximum degree,
		the following upper bound is useful.
		
		\begin{proposition}
		\label{prop:bip2}
			If $G$ is a bipartite graph with  parts  $X$ and $Y$, then
			$$\thetaint(G) \leq \min\{\Delta(X), \Delta(Y)\}.$$
		\end{proposition}
		\begin{proof}
			Let $G$ be a bipartite graph with  parts  $X$ and $Y$.
			We prove that $\thetaint(G) \leq \Delta(X)$; the proof is by induction
			on $\Delta(X)$.
			
			If $\Delta(X)= 1$, then $G$ is a disjoint union of stars, and thus
			$\thetaint(G) =1$.
			If $\Delta(X) \geq k$, then we form a subgraph $H$ of $G$
			by picking one edge incident to each vertex of degree $\Delta(X)$
			in $G$. The graph $H$ is a disjoint union of stars, and hence
			interval colorable. Since $\Delta(G-E(H)) = \Delta(G)-1$,
			the desired result now follows by induction.
		\end{proof}
		
		The upper bound in Proposition \ref{prop:bip2} is tight,
		since there are bipartite graphs without interval colorings
		where all vertices in one part have degree two (see e.g. \cite{PetrosKhacha}).

\subsubsection{An application in  timetabling}
\label{timetable}
In this section we consider a detailed 
application of multi-stage no-wait schedules
 in school timetabling.

In a school,  there are $m$ teachers $P_1,\dots,P_m$ and  $n$
classes $J_1,\dots,J_n$.
A class  consists of a  set of students who follow
exactly the same program.
We are given an $n\times m$  requirement matrix $B=(b_{ij})$ where $b_{ij}$
is the number of
lectures involving class $J_i$ and teacher $P_j$. We shall assume that all
lectures have the same duration (say one period or one hour) 
and for every day of the week the lessons can take place in periods $1,2,3,\dots$. 

A {\em weekly timetable} for $k$ days, corresponding to the requirement matrix $B$, is a sequence  ${\bf S}=(S_1,S_2,\dots,S_k)$ of $k$ arrays with $n$ rows satisfying the following conditions:
\begin{itemize}

\item[(i)]
each entry of $S_l$ is either one of the
members of the set $\{P_1,\dots,P_m\}$ or is  empty, $l=1,\dots,k$;
\item[(ii)]
the total number occurences of the symbol $P_j$ in the $i$th row of the arrays $S_1,\dots,S_k$ is precisely $b_{ij}$,  for $j=1,\dots,m$;

\item[(iii)]
in each column of $ S_l$ all non-empty
symbols are different, $l=1,\dots,k$.
\end{itemize}

Each of the arrays $S_1,\dots,S_k$ is called a {\em daily timetable}. 
If $S_l=(s_{ih}^l)$ and $s_{ih}^l$ is empty then the class $J_i$ has a free lesson in the $h$th period. 
If however $s_{ih}^l=P_j$ then class $J_i$ has a lesson with teacher 
$P_j$ that period.
In fact, if $k=1$ then the array ${\bf S}$ is a daily and weekly timetable simultaneously.
 
Some results on weekly timetables were found byde Werra \cite{deWerra,deWerraD}. In particular he  showed that for any $k\geq 2$ there exists a weekly timetable for 
$k$ days in which the lessons for  each class  and each teacher  are spread throughout the week as evenly as possible.

Let us now %Now we will 
consider timetables without interruptions. 
We say that the class $J_i$ (the teacher $P_j$) has an {\em interruption}
% (a class-interruption)
in a weekly timetable
if there are two periods $j_1$ and $j_2$ at some day of the week, 
such that $j_1+1<j_2$, the class $J_i$  (the teacher $P_j$) 
has lessons at the periods  $j_1$ and $j_2$, but it is 
free at the period $j_1+1$.
  
  If there are no interruption for each of classes and teachers in 
	the daily timetable $S_i$, for  $i=1,\dots,k$, 
	then we say that %${\bf S}=
	$(S_1,S_2,\dots,S_k)$ is a 
	{\em weekly timetable without interruptions}.
   
 Here we are specifically interested in whether there
does exist a weekly timetable 
$(S_1,S_2,\dots,S_k)$   without interruptions, corresponding to the requirement matrix $B$.
The next result shows that this problem is in fact equivalent to Problem \ref{prob:decomp}, stated in the introduction, for a bipartite graph $G=G(B)$ with parts $V_1$ and $V_2$, where
 $V_1=\{J_1,\dots, J_n\}$ and $V_2=\{P_1, \dots,P_m\}$, 
and where the vertices $J_i$ and $P_j$ are joined by $b_{ij}$  edges.

\begin{proposition}
\label{prop:application}
A weekly school timetable for $k$ days without interruptions, corresponding to the requirement matrix $B$,
exists if and only if $\theta_{int}(G(B))\leq k$. 
\end{proposition}

\begin{proof}
 Suppose that there exists a weekly timetable  ${\bf S}=(S_1,S_2,\dots,S_k)$ without interruptions, corresponding to the requirement matrix $B$
 where $S_l=(s_{ih}^l)$, $l=1,\dots,k$. 
Then we define an edge colored subgraph $G_l$ of $G$ as the 
subgraph formed by the edges with ends $J_i$ and $P_j$ colored $h$ 
 under condition that $s_{ih}^l=P_j$, $h=1,2,\dots$. 
 Clearly, such a coloring of $G_l$ is an interval coloring with colors $1,2,3,\dots$.
 
Conversely, 
suppose that the graph $G=G(B)$ can be decomposed into a union of $k$
interval colorable subgraphs $G_1,G_2,\dots,G_k$.
Then we can define a weekly timetable ${\bf S}=(S_1,S_2,\dots,S_k)$ 
with $S_l=(s_{ih}^l)$, for $l=1,\dots,k$, as follows: 
$s_{ih}^l=P_j$ if and only if an edge with ends $J_i$ and $P_j$ in $G_l$ is colored $h$, $l=1,\dots,k$.
\end{proof}

Proposition \ref{th:Bip} can be reformulated as a result 
in terms of timetabling as follows:

\begin{proposition}
\label{prop:even}
Let $B=(b_{ij})$ be an $n \times m$ requirement matrix, let 
$$\Delta = \max \left\{\max_{1\leq j\leq m}\sum_{i=1}^nb_{ij}, 
\max_{1\leq i\leq n}\sum_{j=1}^mb_{ij} \right\},$$
and set $k=\left\lceil \frac{\Delta}{3}\right\rceil$.
Then there is a weekly school timetable for $k$ days without interruptions  
in which the lessons for each class and each teacher are spread throughout the week as evenly as possible.
\end{proposition}

\subsubsection{Biregular graphs}

Now  we consider biregular graphs; a bipartite graph
is {\em $(a,b)$-biregular} if all vertices in one part have degree $a$ and all
vertices in the other part have degree $b$.
Our investigation of biregular graphs is partially motivated by a well-known conjecture
which suggests that all biregular graphs have interval colorings 
\cite{JensenToft}. 
Moreover, biregular graphs arise naturally in some scheduling problems;
for instance, in a school timetable (see section \ref{timetable}) where all teachers have the same number of lectures,
and this also holds for all classes, the problem
of constructing a (weekly) timetable
can be formulated in terms of
an edge coloring problem of a biregular graph.

It is known that all 
$(2,b)$-biregular \cite{Hansen,HansonLotenToft,KamMir} and
$(3,6)$-biregular graphs \cite{CarlJToft} admit interval colorings.
In this section we obtain some upper bounds on the interval coloring thickness
of some biregular graphs.
We shall need the following corollary of Lemma \ref{Lem:add}.

\begin{corollary}\label{Lemma(2,2r)}(\cite{Hansen})
If $G$ is a $(2,2r)$-biregular ($r\geq 2$) bipartite graph 
with parts $X$ and $Y$, then $G$ has an interval 
$2r$-coloring $\alpha$ such that for each $y\in Y$, 
$S_{G}(y,\alpha)=[1,2r]$ and for each
$x\in X$, $S_{G}(x,\alpha)=\{2i-1,2i\}$ for some $i$,
$1\leq i\leq r$.
\end{corollary}

\begin{proposition}
\label{prop:bireg(3,3r)and(4,4r)}
Let $G$ be a bipartite graph.
\begin{itemize}
	
	\item[(i)] If $G$ is $(3,3r)$-biregular ($r\geq 2$), 
	then $\theta_{int}(G)\leq 2$.

	\item[(ii)] If $G$ is $(k,kr)$-biregular ($k \geq 4, r\geq 2$),
	then $\theta_{int}(G)\leq k-2$.
	
\end{itemize}
\end{proposition}
\begin{proof}
We first prove the upper bound in (i).  

Let $G$ be a $(3,3r)$-biregular bipartite graph with parts $X$ and $Y$.
Define a new graph $H$ from $G$ by splitting each vertex $y\in Y$ into 
$r$ vertices
$y^{(1)},y^{(2)},\ldots,y^{(r)}$ of degree $3$.
The graph $H$ is cubic and bipartite, so by Hall's matching theorem, 
it has a perfect matching $M$.

In the graph $G$, $M$ induces a subgraph $F$ in which each vertex $y\in Y$ has degree $r$ and each vertex $x\in X$ has degree $1$, so $F$ is interval colorable.
Moreover,
since $G^{\prime} = G-E(F)$ is $(2,2r)$-biregular, it follows from 
Corollary \ref{Lemma(2,2r)} that it has an interval $2r$-coloring. 
Thus, $\theta_{int}(G)\leq 2$.

Next, we prove the upper bound in (ii).
The proof is by induction on $k$. Let us first consider the
base case when $k=4$.

Without loss of generality, we may assume that $G$ is a connected
(otherwise, we consider every connected component of $G$)
$(4,4r)$-biregular graph. Since $G$
is bipartite and all vertex degrees in $G$ are even, $G$ has a
closed Eulerian trail $C$ with an even number of edges. We color the
edges of $G$ with colors red and blue by traversing the edges
of $G$ along the trail $C$; we color an odd-indexed edge in $C$ with
color red, and an even-indexed edge in $C$ with color blue. 
Let $G_R$ be the subgraph induced by the red edges, and
$G_B$ be the graph induced by the blue edges.

Since $G$ is $(4,4r)$-biregular, each of the subgraphs $G_{R}$ and $G_{B}$
is a $(2,2r)$-biregular graph. Thus,
by Corollary \ref{Lemma(2,2r)}, $\theta_{int}(G)\leq 2$.

Let us now assume that $k\geq 5$ and that the statement is true 
for any $(k',k'r)$-biregular bipartite graph $G'$, where $k'<k$.
Let $G$ be a $(k,kr)$-biregular bipartite graph with parts $X$  and $Y$.
We define a new graph $H$ from $G$ by splitting each vertex $y\in Y$ into
$r$ vertices $y^{(1)},y^{(2)},\ldots,y^{(r)}$ of degree $k$.
The graph $H$ is a $k$-regular and bipartite, so by Hall's 
matching theorem, $H$ 
contains a perfect matching $M$. 

In the graph $G$, $M$ induces a subgraph $F$ in which each 
vertex of $Y$ has degree $r$ and each vertex of $X$ has degree $1$,
so $F$ is interval colorable. Moreover, 
the graph $G^{\prime} = G-E(F)$ is $((k-1),(k-1)r)$-biregular,
so by the induction hypothesis, $\theta_{int}(G^{\prime})\leq k-3$. 
Thus, $\theta_{int}(G)\leq k-2$.
\end{proof}

\subsection{Decompositions into forests and complete multipartite graphs}

The {\em arboricity} of a graph $G$, denoted $\gamma(G)$, 
is the least number of edge-disjoint  forests whose union is $G$.
Since all
trees and forests are interval colorable, $\thetaint(G)\leq \gamma(G)$ 
for every graph $G$. Thus
the  arboricity $\gamma(G)$ immediately yields an upper bound on $\thetaint(G)$.
It was proved by Nash-Williams \cite{NashWilliams} that
$$\gamma(G)=\max \left\lceil\frac{|E(G[X])|}{|X|-1}\right\rceil$$
where the maximum is taken over all nonempty subsets $X\subseteq V(G)$.
This result implies several different upper bounds on $\gamma(G)$. 
We note below only two of them.

\begin{theorem} (\cite{DeanHutchinsonScheinerman})
\label{th:Dean}
	If $G$ is a simple graph with $q$ edges then $\gamma(G)\leq  
	\left\lceil\sqrt {q/2}\right\rceil$. 
	\end{theorem}

 The following is well-known.

\begin{proposition} 
\label{th:plane}
	If $G$ is a simple planar  graph  then $\gamma(G) \leq 3$.   
	If, additionally, $G$ is triangle-free 
	or outerplanar, then 
	$\gamma(G) \leq 2$.
\end{proposition}

Theorem \ref{th:Dean} and Proposition \ref{th:plane} together with the
fact that $\thetaint(G)\leq \gamma(G)$ imply the following for simple graphs:

\begin{proposition}
\label{prop:dean}
\begin{itemize}

	\item[(i)] 
		If $G$ is a simple graph with $q$ edges, then 
		$\theta_{int}(G)\leq  \left\lceil\sqrt {q/2}\right\rceil$. 
		
	\item[(ii)] If $G$ is a simple planar  graph  then $\thetaint(G) \leq 3$.  
	If, additionally, $G$ is triangle-free 
	or outerplanar, then 
	$\thetaint(G) \leq 2$.
	
		\end{itemize}
	\end{proposition}
	
	Note that this upper bound for outerplanar graphs is sharp, since
	odd cycles are not interval colorable, while it is an open question
	whether there are planar graphs with interval coloring thickness $3$.

Note that for some graphs $G$ the number $\gamma(G)-\thetaint(G)$
can be very large. Consider, for example, the complete 
graph $K_{2n+1}$, $n\geq 1$. It is known \cite{AsrKam} that $\thetaint(K_{2n+1}) > 1$. 
Moreover, for any 
$u \in V(G)$, 
		the graph $H=K_{2n+1}-u$ is a complete graph with $2n$ vertices
		which is interval colorable. Since $K_{2n+1} - E(H)$ is a tree, which is an interval colorable graph, 
		 $\thetaint(K_{2n+1})=2$.
		
On the other hand, by the formula of Nash-Williams for $\gamma(G)$ above 
$$\gamma(K_{2n+1})\geq \left\lceil \frac{|E(K_{2n+1})|}{|V(K_{2n+1})|-1}\right\rceil\geq \left\lceil \frac{(2n+1)n}{2n}\right\rceil=
\left\lceil n+\frac{1}{2}\right\rceil=n+1.$$   
		
		Therefore, for any $n \geq 1$, $\gamma(K_{2n+1})-\theta(K_{2n+1})\geq n-1$.

	\bigskip

	Instead of decomposing a graph $G$ into trees, we could consider 
	decompositions $G$ into
	complete bipartite graphs, which are also known to always 
	admit interval colorings.
	We shall demonstrate this method for the case of complete multipartite graphs.

	A simple graph is called {\em complete $r$-partite} ($r\geq 2$) if its vertices can be partitioned into $r$ nonempty independent sets $V_1,\ldots,V_r$ such that each vertex in $V_i$ is adjacent to all the other vertices in $V_j$ for $1\leq i<j\leq r$. Let $K_{n_{1},n_{2},\ldots,n_{r}}$ denote a complete $r$-partite graph
with independent sets $V_1,V_2,\ldots,V_r$ of sizes $n_{1},n_{2},\ldots,n_{r}$.

Using the following proposition, 
which follows from the so-called Master
theorem (see e.g. \cite{Cormen}), we shall deduce an upper 
bound on $\thetaint(K_{n_{1},n_{2},\ldots,n_{r}})$ for arbitrary $n_1,\dots, n_r$.

\begin{proposition}
 \label{cor:recurrence}
 If $T(n)$ is a function defined by the recurrence
$$T(n)= T(\left\lceil{n/2}\right\rceil)+1,$$
then $T(n)=\Theta(\log n)$. 
 \end{proposition}
	 	
First, we note the following proposition, which follows from the fact that
a complete $4$-partite graph can be decomposed into two subgraphs, 
each consisting of two
disjoint complete bipartite graphs.

\begin{proposition}
\label{prop:4part}
	If $G$ is a complete multipartite graph with at most
	four parts, then $\thetaint(G) \leq 2$.
\end{proposition}

Note that Proposition \ref{prop:4part} is sharp since there are complete
$4$-partite graphs of Class 2.

\begin{theorem}
\label{completemulti} 
 If $G$ is a complete $r$-partite graph, then $\theta_{int}(G)= O(\log r)$.
\end{theorem} 

\begin{proof}
Let  $V_1,\ldots,V_r$ be nonempty independent sets in $G=K_{n_{1},n_{2},\ldots,n_{r}}$ such that each vertex in $V_i$ is adjacent to all the other vertices in $V_j$ for $1\leq i<j\leq r$. 
For any set $A\subset V(G)$ that is a union of some $t\geq 2$ 
subsets $V_{i_1},\dots,V_{i_t}$ we  denote by $F(A)$ the complete bipartite graph with parts 
 $V_{i_1}\cup...\cup V_{i_k}$ and $V_{i_{k+1}}\cup...\cup V_{i_t}$ 
where $k=\left\lceil\frac{t}{2}\right\rceil$. 

By Proposition \ref{cor:recurrence}, we may assume that
$r\geq 5$. We construct a decomposition  of $G$ 
into edge-disjoint graphs $Q_1, Q_2,\dots$  where the components of every $Q_i$ 
are disjoint complete
bipartite graphs as follows:

Step 1 . Set $Q_1=F(V_1\cup V_2\cup\dots\cup V_r)$. 

Note that the subgraph  of $G$ induced by the set $V_1\cup\dots\cup V_k$ 
is a complete $k$-partite graph and the subgraph induced by the set $V_{1+k}\cup\dots\cup V_r$ is a complete $(r-k)$-partite graph, where $k=\left\lceil\frac{r}{2}\right\rceil$.

Step 2. Set $Q_2=F(V_1\cup\dots\cup V_k)\cup F(V_{1+k}\cup\dots\cup V_r)$ where $k=\left\lceil\frac{r}{2}\right\rceil$.  

Step $i$ ($i\geq 3$). Suppose that the graphs $Q_1,\dots,Q_{i-1}$ 
 have been constructed and $Q_{i-1}=H_1\cup\dots\cup H_l$  where $H_j$ is
a complete bipartite graph with parts  $A_j$  and $B_j$,
$j=1,\dots,l$, such that all sets $A_1,\dots, A_l$, $B_1,\dots, B_l$ are mutually disjoint. 

 If every $A_j$ is one of the sets $V_1,\dots,V_r$
and every $B_j$ is also one of the sets $V_1,\dots,V_r$, $j=1,\dots,l$, then STOP the algorithm. 
Otherwise we construct the bipartite graph $F(A_j)$ for each set $A_j$  containing at least 2 of subsets $V_1,\ldots,V_r$, $1\leq j\leq l$,
and  the bipartite graph $F(B_s)$ for each set $B_s$  containing at least 2 of subsets $V_1,\ldots,V_r$, $1\leq s\leq l$.
(Note that each such set induces a complete multipartite graphs.) 
Now define $Q_i$ as a union of all these bipartite graphs. Formally, we set
$$Q_i=\left(\cup_{A_j\notin \{V_1,\dots,V_r\}}F(A_j)\right)\cup 
\left(\cup_{B_s\notin \{V_1,\dots,V_r\}}F(B_s)\right).$$

This algorithm stops when  all edges of $G$ will be 
 in one of the complete bipartite graphs constructed by the algorithm.
We denote by $T(r)$ the number of bipartite graphs $Q_1,Q_2,\dots$ constructed by this algorithm. Clearly, 
$$T(r)=T(\left\lceil\frac{r}{2}\right\rceil)+1,$$
so by  Proposition \ref{cor:recurrence}, $T(r)=\Theta (\log r)$.

Since every component of the graph $Q_i$ in the constructed decomposition of $G=K_{n_{1},n_{2},\ldots,n_{r}}$ is a complete bipartite graph, all $Q_i$ are
 interval colorable,
and
therefore $\theta_{int}(G)\leq T(r)=\Theta (\log r)$. 
\end{proof}

Finally,  let us note some further upper bounds on the interval coloring thickness
of complete multipartite graphs.
A complete $r$-partite graph
$K_{n_{1},n_{2},\ldots,n_{r}}$ is called {\em balanced} and  denoted by
  $K_{n*r}$, if $n_1=\dots=n_r=n$.

Also, let $K_{n*r,a}$ denote the complete $(r+1)$-partite graph with $n$ vertices in each of the first $r$ parts and $a$ vertices in the last part.
In \cite{PetArXiv}, the following consequence on the chromatic index of a
balanced complete multipartite graph was observed.

\begin{proposition}
\label{complete-regular-int}
For any $n,r\in \mathbb{N}$ ($r\geq 2$), $K_{n*r}$ has an interval coloring if and only if $nr$ is even. Moreover, if $nr$ is even, then $K_{n*r}$ has an interval $(r-1)n$-coloring. 
\end{proposition}

Using the preceding proposition, we deduce the following two results.

\begin{proposition}
\label{complete-regular} 
For any $n,r\in \mathbb{N}$ ($r\geq 2$), 
\begin{center}
$\theta_{int}\left(K_{n*r}\right)=\left\{
\begin{tabular}{ll}
$1$, & if $nr$ is even,\\
$2$, & if $nr$ is odd.\\
\end{tabular}%
\right.$
\end{center}
\end{proposition}

Proposition \ref{complete-regular} follows from Proposition \ref{complete-regular-int}
together with the fact that if $nr$ is odd then, $K_{n*r}$ 
can be decomposed into a copy of $K_{n*(r-1)}$ and a complete bipartite graph.

\begin{proposition}
\label{complete-semiregular} 
For any $n,r\in \mathbb{N}$ ($r\geq 2$),
$\theta_{int}\left(K_{n*r,nr}\right)\leq 3$. Moreover, if $nr$ is even, then $\theta_{int}\left(K_{n*r,nr}\right)=1$.
\end{proposition} 
The upper bound follows from the previous proposition, since
$K_{n*r,nr}$ can be decomposed into a copy of $K_{n*r}$
and a balanced complete bipartite graph.

The second part follows by considering the same decomposition
and using Proposition \ref{complete-regular} for obtaining an interval $(r-1)n$-coloring
of the copy of  $K_{n*r}$, and then defining a suitable interval coloring of
the balanced complete bipartite graph.

\section{Concluding remarks and open problems}

In this paper we have introduced and investigated the interval 
coloring thickness of a graph
and considered some concrete applications from scheduling theory. We have presented
a new class of interval colorable graphs and proved an upper
bound on the interval coloring
thickness of a general graph. Furthermore, we have proved
improvements of this general upper
bound for several families of graphs.
Let us now point to some open problems,
and also describe
some further tractable instances.

There are several well-known examples of families of graphs that
do not admit interval colorings, such as odd cycles
and complete graphs of odd order.
However, the answer to the following question
remains open.

\begin{problem}
\label{prob:general}
	For any positive integer $k$, is there a
	graph
	$G$ with $\theta_{int}(G) =k$?
\end{problem}

The version of this problem when the graph is assumed to be bipartite is
perhaps of particular interest.
In \cite{GiaroKubaleMalaf1,PetrosKhacha,AsratianCasselgrenPetrosyan1}, the authors consider a variety of families of bipartite graphs
without interval colorings; 
these families include, the the so-called Malafiejski rosettes, the Sevastjanov rosettes
and the generalized Hertz graphs. All graphs from these families have interval coloring
thickness two, since they all decompose into two forests.
Note further that any example
of a bipartite graph with interval coloring thickness $k$ must necessarily
contain at least $3k+1$ vertices in each part, since every bipartite graph
where one part contains at most three vertices is interval colorable
\cite{GiaroKubaleMalaf1}.

Based on an example first described by Erd\"os, in \cite{PetrosKhacha} the authors present
a family of bipartite graphs, constructed from projective planes, that do not admit interval
colorings. The interval coloring thickness of graphs from this family is unknown.
A solution to the following problem could perhaps shed some light on Problem \ref{prob:general}.

\begin{problem}
\label{Erd}
	Determine the interval coloring thickness of the so-called Erd\"os family of graphs
	described in Section 2.2. in \cite{PetrosKhacha}.
\end{problem}

Some  graphs which are not interval colorable, permit  
a proper edge coloring that satisfy a slightly weaker condition.
A proper $t$-edge coloring of a graph $G$ is
%$\alpha:E(G)\longrightarrow  \{1,\dots,t\}$ of a graph $G$  is
called a \emph{cyclic interval $t$-coloring} if for each vertex $v$ of $G$ the edges incident to $v$ are colored by consecutive colors, under the condition that color $1$ is considered as consecutive to color $t$, 
see e.g. \cite{AsratianCasselgrenPetrosyan}.
We note the following proposition for graphs with a cyclic interval coloring.

\begin{proposition}
\label{cyclicthm}
If $G$ has a cyclic interval $t$-coloring with $t\geq 2\Delta(G)-2$, then 
$\theta_{int}G)\leq 2$.
\end{proposition}

Proposition \ref{cyclicthm} can be proved as follows:
Let $\alpha$ be a cyclic interval $t$-coloring of $G$; 
if $\alpha$ is an interval coloring, then we are done;
otherwise, there is some vertex
 $v\in V(G)$
such that $S(v,\alpha)$ is not an interval, that
is $S(v,\alpha)=\{1,\ldots,k_v\}\cup
\{t-l_v+1,\ldots,t\}$, for some integers $k_v,l_v$.
Let $k^*$ be the maximum integer $k_x$ such
that $S(x,\alpha)=\{1,\ldots,k_x\}\cup
\{t-l_x+1,\ldots,t\}$, for some vertex $x$,
and let $l^*$ be the maximum integer $l_y$ such that
$S(y,\alpha)=\{1,\ldots,k_y\}\cup
\{t-l_y+1,\ldots,t\}$, for some vertex $y$ (where $S(x,\alpha)$ and $S(y,\alpha)$ are 
not intervals).
Since $t \geq 2 \Delta(G)-2$, $k^* < t-l^*+1$, and thus
$G$ decomposes into two interval colorable subgraphs.

An example of a family of graphs satisfying Proposition \ref{cyclicthm}
is the so-called {\em circular complete graphs}
$K_{p/q}$ defined by setting
$V(K_{p/q})=\{v_{0},v_{1},\ldots,v_{p-1}\}$ and $$E(K_{p/q})=\{v_{i}v_{j}: v_{i},v_{j}\in V(K_{p/q}) \text{~and~} q\leq |i-j|\leq p-q\}.$$ 
In general, we are interested in whether the following might be true.

\begin{problem}
\label{cyclic}
	Is it true that every graph $G$ that admits a cyclic interval coloring satisfies 
	$\thetaint(G) \leq 2$?
\end{problem}

\end{document}